\documentclass[12pt]{amsart}
\usepackage{natbib}

\theoremstyle{definition}
\newtheorem{definition}{Definition}[section]

\theoremstyle{plain}
\newtheorem{theorem}[definition]{Theorem}
\newtheorem{lemma}[definition]{Lemma}

\newtheorem{proposition}[definition]{Proposition}

\theoremstyle{remark}

\def\to{\longrightarrow}
\def\paige#1{M^*(#1)}
\def\vm#1#2#3#4{\begin{pmatrix} #1 & #2 \\ #3 & #4 \end{pmatrix}}
\def\psl#1#2{L_#1(#2)}
\def\sl#1#2{SL_#1(#2)}
\def\octo{\mathbb O}
\def\field#1{GF(#1)}
\def\dpr#1#2{#1\cdot#2}
\def\vpr#1#2{#1\times#2}
\def\centre#1{Z(#1)}

\begin{document}


\title[Generators for finite simple Moufang loops]
{Generators for Finite Simple Moufang Loops
    }
\author{Petr Vojt\v echovsk\'y}
\address{Department of Mathematics, Iowa State University, Ames, IA, U.S.A.}
\email{petr@iastate.edu}

\begin{abstract}
Moufang loops are one of the best-known generalizations of groups. There is
only one countable family of nonassociative finite simple Moufang loops,
arising from the split octonion algebras. We prove that every member of this
family is generated by three elements, using the classical results on
generators of unimodular groups.

\end{abstract}

\maketitle

\section{Finite Simple Moufang Loops}

\noindent A \emph{Moufang loop} is a loop satisfying the identity
$x(y(xz))=((xy)x)z$. Every element $x$ of a Moufang loop has a (unique)
two-sided inverse $x^{-1}$; every two elements generate a group. Loops with the
latter property are called \emph{diassociative} \cite{Pflugfelder}. The
best-known Moufang loop is the multiplicative loop of nonzero elements in the
standard $8$-dimensional real octonion algebra $\octo$. Surely the best-known
\emph{finite} Moufang loop is the $240$-element loop $L$ of integral octonions
of norm one \cite{Coxeter}.

In $1956$, L.~Paige \cite{Paige} found one nonassociative finite simple Moufang
loop for every finite field. Following E.~Bannai and S.~Song \cite{BannaiSong},
we denote this \emph{Paige loop} constructed over $F=\field{q}$ by $\paige{q}$.
Let us give the most brief description of $\paige{q}$ now.

Consider the Zorn multiplication
\begin{equation}\label{Eq:Zorn}
    \vm{a}{\alpha}{\beta}{b}\vm{c}{\gamma}{\delta}{d}
        =\vm{ac+\dpr{\alpha}{\delta}}{a\gamma+\alpha d-\vpr{\beta}{\delta}}
            {\beta c+b\delta+\vpr{\alpha}{\gamma}}{\dpr{\beta}{\gamma}+bd},
\end{equation}
where $a$, $b$, $c$, $d\in F$, $\alpha$, $\beta$, $\gamma$, $\delta\in
F^{\,3}$, and where $\dpr{\alpha}{\delta}$ (resp. $\vpr{\alpha}{\delta}$)
denotes the dot product (resp. cross product) of $\alpha$ and $\delta$. This is
the same formula M.~Zorn used to construct the split octonion algebra over $F$.
The loop $\paige{q}$ consists of all matrices
\begin{displaymath}
    M = \vm{a}{\alpha}{\beta}{b}
\end{displaymath}
with $\det{M}=ab-\alpha\beta=1$ that are multiplied according to
$(\ref{Eq:Zorn})$, and where $M$ and $-M$ are identified. The neutral element
of $\paige{q}$ is the identity matrix $I$, and the inverse of $M$ is
\begin{displaymath}
    M^{-1}=\vm{b}{-\alpha}{-\beta}{a}.
\end{displaymath}

In $1987$, M.~Liebeck \cite{Liebeck} proved that there are no other
nonassociative finite simple Moufang loops. The loop $\paige{2}$ is exceptional
in the sense that it shows up in the real algebra $\octo$, too. Namely,
$\paige{2}$ is isomorphic to the quotient of $L$ by its center
$\centre{L}=\{1,\,-1\}$ (see \cite{Coxeter} once again).

Associative finite simple Moufang loops are finite simple groups. It is a
remarkable fact that every finite simple group is $2$-generated
\cite{AschbacherGuralnick}; even more so, since no proof using only the
simplicity is known. Instead, every family of finite simple groups must be
investigated separately. Because of diassociativity, the nonassociative Paige
loops cannot be $2$-generated. It is reasonable to expect that a small number
of generators will do. Indeed, it this paper we prove that:

\begin{theorem}\label{Th:Main}
Every nonassociative finite simple Moufang loop is $3$-generated.
\end{theorem}

Note that Theorem \ref{Th:Main} was proved in \cite{JA} for all Paige loops
$\paige{p}$, $p$ a prime. Thus the main task of this paper is to cover the
general case. We also present a simple proof for the prime case, and offer at
least two generating sets for every $\paige{q}$. The reader who wants to
establish \ref{Th:Main} as quickly as possible should focus on
$\ref{Pr:AlbertThompson}$, $\ref{Lm:CoxeterMoser}$, and Section
$\ref{Sc:Paige}$.

\section{Generators for $\psl{2}{q}$}

\noindent The crucial observation concerning Paige loops is that $\paige{q}$
contains several copies of $\psl{2}{q}=PSL_2(q)$. Given the canonical basis
$e_1=(1,\,0,\,0)$, $e_2=(0,\,1,\,0)$, $e_3=(0,\,0,\,1)$ of $F^{\,3}$, let
$\phi_i:\psl{2}{q}\to\paige{q}$ be defined by
\begin{displaymath}
    \phi_i\vm{a}{b}{c}{d}=\vm{a}{be_i}{ce_i}{d},
\end{displaymath}
and let $G_i$ be the image of $\psl{2}{q}$ under $\phi_i$. Since the
multiplication in $G_i$ coincides with the usual matrix multiplication (all
cross products involved in $(\ref{Eq:Zorn})$ vanish), $\phi_i$ is an
isomorphism $\psl{2}{q}\to G_i$.

This brings our attention to the classical results concerning generators for
$\psl{2}{q}$ and $\sl{2}{q}$. First of all, we have the Dickson Theorem:

\begin{theorem}[L.~E.~Dickson, 1900]\label{Th:Dickson}
If $q\ne 9$ is an odd prime power or $q=2$, then $\sl{2}{q}$ is generated by
\begin{equation}\label{Eq:Dickson}
    \vm{1}{1}{0}{1},\;\;\vm{1}{0}{\lambda}{1},
\end{equation}
where $\lambda$ is a primitive element of $\field{q}$.
\end{theorem}

The proof can be found in \cite{Dickson}, and more recently in \cite[pp.\
$44$--$55$]{Gorenstein}. The statement of the theorem usually does not mention
$q=2$, although it is apparently true for $q=2$, since $\sl{2}{2}\cong S_3$ is
generated by any two involutions, in particular by $(\ref{Eq:Dickson})$.

A.~A.~Albert and J.~Thompson proved \cite[Lemma $8$]{AlbertThompson} that for
any primitive element $\lambda$ of $\field{q}$, $q>2$, the group $\sl{2}{q}$ is
generated by $B$, $-B$, and $C$, where
\begin{equation}\label{Eq:AlbertThompson}
    B=\vm{\lambda}{0}{0}{\lambda^{-1}},\;\; C=\vm{0}{1}{-1}{\lambda}.
\end{equation}
We therefore have:

\begin{proposition}[A.~A.~Albert, J.~Thompson, 1959]\label{Pr:AlbertThompson}
Let $q>2$ be a prime power. Then $\psl{2}{q}$ is generated by
$(\ref{Eq:AlbertThompson})$, where $\lambda$ is a primitive element of
$\field{q}$.
\end{proposition}

The generators $(\ref{Eq:AlbertThompson})$ are especially convenient for our
purposes, because $\phi_i(B)=B$ for every $i$, $1\le i\le 3$; but let us not
get ahead of ourselves. It is practical to know some generators that do not
involve a primitive element. For that matter, Coxeter and Moser argue in
\cite{CoxeterMoser} that

\begin{lemma}\label{Lm:CoxeterMoser}
For every prime $p$, the group $\psl{2}{p}$ is generated by
\begin{equation}\label{Eq:CoxeterMoser}
    \vm{1}{0}{1}{1},\;\;\vm{0}{1}{-1}{0}.
\end{equation}
\end{lemma}

\section{Generators for $\paige{q}$}\label{Sc:Paige}

\noindent Our first result concerning $\paige{q}$ has nothing to do with the
generators for $\psl{2}{q}$. In its proof, we take advantage of the following
lemma due to Paige:

\begin{lemma}[L.~Paige, 1956]\label{Lm:Paige}
$\paige{q}$ is generated by
\begin{equation}\label{Eq:PaigeGens}
    M_\beta=\vm{1}{0}{\beta}{1},\;\; M'_\beta=\vm{1}{\beta}{0}{1},
\end{equation}
where $\beta$ runs over all nonzero vectors in $F^{\,3}$.
\end{lemma}
\begin{proof}
Combine Lemmas $4.2$ and $4.3$ of \cite{Paige}.
\end{proof}

\begin{proposition}\label{Pr:G1G2G3}
$\paige{q}$ is generated by $G_1\cup G_2\cup G_3$.
\end{proposition}
\begin{proof}
Let $Q$ be the subloop of $\paige{q}$ generated by $G_1\cup G_2\cup G_3$.
Thanks to Lemma \ref{Lm:Paige}, it suffices to prove that $Q$ contains all
elements $M_\beta$, $M'_\beta$, defined in $(\ref{Eq:PaigeGens})$. We show
simultaneously that $M_\beta\in Q$ and $M'_\beta\in Q$.

Let $k$ denote the number of nonzero entries of $\beta$. There is nothing to
prove when $k\le 1$. Suppose that $k=2$. Without loss of generality, let
$\beta=(a$, $b$, $0)$ for some $a$, $b\in F^*=F\setminus\{0\}$. Verify that
\begin{displaymath}
    \vm{1}{ae_1}{0}{1}\vm{1}{be_2}{0}{1}\cdot\vm{1}{0}{-abe_3}{1}
        =\vm{1}{(a,b,0)}{0}{1},
\end{displaymath}
and thus that $M_\beta\in Q$. Similarly, $M'_\beta\in Q$. We can therefore
assume that $Q$ contains all elements $M_\beta$, $M'_\beta$ with $k\le 2$.

Let $k=3$, $\beta=(a$, $b$, $c)$ for some $a$, $b$, $c\in F^*$. As
\begin{displaymath}
    \vm{1}{(a,b,0)}{0}{1}\vm{1}{(0,0,c)}{0}{1}\cdot\vm{1}{0}{(-bc,ac,0)}{1}
        =\vm{1}{(a,b,c)}{0}{1},
\end{displaymath}
$M_\beta$ belongs to $Q$. Symmetrically, $M'_\beta\in Q$, and we are done.
\end{proof}

In fact, $G_1\cup G_2$ already generates $\paige{q}$. The role of the cross
product is especially apparent in the next Proposition.

\begin{proposition}\label{Pr:G1G2}
The subgroup $G_3$ is contained in the subloop of $\paige{q}$ generated by
$G_1\cup G_2$. In particular, $\paige{q}$ is generated by $G_1\cup G_2$.
\begin{proof}
As it turns out, all we need are these two equations:
\begin{eqnarray*}
    \vm{1}{0}{\lambda e_3}{1}&=&
        -\vm{0}{e_2}{-e_2}{0}\vm{1}{\lambda e_1}{-\lambda^{-1}e_1}{0}
        \cdot\vm{1}{e_2}{-e_2}{0}\vm{1}{\lambda e_1}{-\lambda^{-1}e_1}{0},\\
    \vm{0}{e_3}{-e_3}{0}&=&
        \vm{0}{e_1}{-e_1}{0}\vm{0}{-e_2}{e_2}{0}.
\end{eqnarray*}
Note that the left hand sides of these equations are elements of $G_3$, whereas
the right hand sides are products of elements of $G_1\cup G_2$. When $q=2$, we
are done by Lemma \ref{Lm:CoxeterMoser}. When $q>2$, observe that
\begin{displaymath}
    \vm{1}{0}{\lambda}{1}\vm{0}{1}{-1}{0}=\vm{0}{1}{-1}{\lambda}=C.
\end{displaymath}
Since $B=\phi_i(B)$ for every $i$, $1\le i\le 3$, we are done by Proposition
\ref{Pr:AlbertThompson}.
\end{proof}
\end{proposition}

Theorem \ref{Th:Main} is now proved. When $q>2$, $\paige{q}$ is generated by
$\phi_1(C)$, $\phi_2(C)$ and $B=\phi_1(B)=\phi_2(B)$, by Propositions
\ref{Pr:AlbertThompson} and \ref{Pr:G1G2}. When $q=2$, we are done by the main
result of \cite{JA}, Theorem $2.1$ \cite{JA}.

For the sake of completeness, allow us to present an alternative, simpler proof
of \cite[Theorem $2.1$]{JA}.

\begin{proposition}\label{Pr:p}
\emph{\cite[Theorem $2.1$]{JA}} Let $p$ be a prime. Then $\paige{p}$ is
generated by
\begin{displaymath}
    U_1=\vm{1}{e_1}{0}{1},\;\;U_2=\vm{1}{e_2}{0}{1},\;\;
    X=\vm{0}{e_3}{-e_3}{1}.
\end{displaymath}
\begin{proof}
First check that
\begin{equation}\label{Eq:Conj}
    \vm{1}{0}{1}{1}=\vm{0}{1}{-1}{0}\vm{1}{1}{0}{1}^{-1}
        \vm{0}{1}{-1}{0}^{-1}.
\end{equation}
Combine $(\ref{Eq:CoxeterMoser})$ and $(\ref{Eq:Conj})$ to see that
$\psl{2}{p}$ is generated by
\begin{displaymath}
    U=\vm{1}{1}{0}{1},\;\;V=\vm{0}{1}{-1}{0}.
\end{displaymath}
Consequently, $\paige{p}$ is generated by $U_1=\phi_1(U)$, $U_2=\phi_2(U)$,
$V_1=\phi_1(V)$, and $V_2=\phi_2(V)$. Now,
\begin{eqnarray*}
     V_2&=&-(XU_1\cdot XU_2)\cdot X^{-1}U_1,\\
     V_1&=&-U_1U_2\cdot(V_2\cdot U_1X),
\end{eqnarray*}
and we are through.
\end{proof}
\end{proposition}

\section{More Generating Sets}

\noindent We would like to show how to obtain additional generating sets for
$\paige{q}$. We take advantage of Proposition \ref{Pr:G1G2}, Dickson's Theorem,
and of the fact that $\sl{2}{2^r}$ (for $r>1$) is generated by
\begin{equation}\label{Eq:2r}
    D_1=\vm{1}{1}{1}{0},\;\;D_2=\vm{\lambda}{0}{0}{\lambda^{-1}},
\end{equation}
where $\lambda$ is a primitive element of $\field{2^r}$. We leave the
verification of $(\ref{Eq:2r})$ to the reader.

Since $\phi_i(D_2)=D_2$ for $i=1$, $2$, $3$, we immediately see from
Proposition \ref{Pr:G1G2} that $\paige{2^r}$ (for $r>1$) is generated by
$\phi_1(D_1)$, $\phi_2(D_1)$ and $D_2$.

\begin{proposition}\label{Pr:Dickson}
Let $q\ne 9$ be an odd prime power or $q=2$. Then $\paige{q}$ is generated by
\begin{equation*}
    \vm{1}{e_1}{0}{1},\;\vm{1}{e_2}{0}{1},\;
        \vm{0}{\lambda e_3}{-\lambda^{-1}e_3}{1},
\end{equation*}
where $\lambda$ is a primitive element of $GF(q)$.
\begin{proof}
Keeping Proposition \ref{Pr:G1G2} and Dickson's Theorem in mind, we only need
to obtain the elements
\begin{displaymath}
    \vm{1}{0}{\lambda e_i}{1},
\end{displaymath}
for $i=1$, $2$. Straightforward computation reveals that
\begin{eqnarray*}
    \vm{1}{0}{\lambda e_1}{1}&=&
        -\vm{0}{\lambda e_3}{-\lambda^{-1}e_3}{1}^2\vm{1}{e_2}{0}{1}
        \vm{0}{\lambda e_3}{-\lambda^{-1}e_3}{1},\\
    \vm{1}{0}{\lambda e_2}{1}^{-1}&=&
        -\vm{0}{\lambda e_3}{-\lambda^{-1}e_3}{1}^2\vm{1}{e_1}{0}{1}
        \vm{0}{\lambda e_3}{-\lambda^{-1}e_3}{1}.
\end{eqnarray*}
Note that the expressions on the right hand side can be evaluated in any order.
\end{proof}
\end{proposition}

\bibliographystyle{plain}

\end{document}